\documentstyle[12pt]{article}
\newtheorem{theorem}{Theorem}

\newtheorem{corollary}{Corollary}

\def\Proof{{\bf Proof}}

\font\Bbb=msbm10 scaled 1200
\def\Z{\hbox{\Bbb Z}}
\def\N{\hbox{\Bbb N}}
\def\T{\hbox{\Bbb T}}
\def\P{\hbox{\Bbb P}}

\def\R{\hbox{\Bbb R}}

\def\beq{\begin{equation}} \def\eeq{\end{equation}}
\font\gothic=eufm10
\def\goth#1{\hbox{\gothic #1}}

\begin{document}

\begin{center}
              QUANTIZED GEODESIC FLOW AND SUPERCONDUCTIVITY OF PLASMA IN FIREBALLS
\end{center}

\bigskip

One looks for superconductivity at room temperature. But the
superconductivity is a purely quantum phenomenon and it is worthwhile
to search it in situations in which the quantum nature of substances
is most clearly defined. These are, naturally, the low temperature,
and besides the state of plasma for comparatively high temperature. I
belive that plasma can be a superconductor and the verification of
this suggestion is the phenomenon of fireball.

In the work "On the nature of fireball" \cite{Kap} P.L.Kapitza
calculated that the whole energy which can be contained in a ball with
the size of a standard fireball (even for the case of full ionization)
would irradiated during the time of order 0.01 sec. But the life time
of fireballs is several minutes, that has increased the theoretical
value ten thousend times. Because of this P.L.Kapitza was forced to
propose the hypothesis of constant feeding fireballs by the energy of
external magnetic field. But the succeeding experiments of generation
such fields by Super-High-Frequency vibrations (SHF) to derive
fireballs have been unsuccessful. There were tens of efforts to find a
physical mechanism which could explain such a long life time of
fireballs. It is generally agreed that all of them collapsed. It seems
to me that there is the only way out of the difficulty. To accept that
the groundwork of the phenomenon is the superconductivity of plasma.

B.B.Kadomzev in the book "Collective phenomena in plasma" \cite{Kad}
considers oscillations in plasma with negative energy which is
indirect evidence for a gap in energy spectrum and hence for the
possibility of superconductivity. In the paper of collaborators of
Lawrens Livermore National Laboratory "Temperature Measurements of
Shock Compreesed Liquid Deuterium up to 230 GPa" \cite{Liv} it was
described that at high pressures and temperatures deuterium behaves
like a degenerate Fermi-liquid metal: its compressibility and
conductivity abruptly increased, it appears a highly reflective state
which is caracteristic of a liquid metal. Such behaviour also points
out to a possibility of superconductivity. An added reason for
superconductivity is the Meisner-Ochsenfeld effect. The
fact that magnetic fields penetrate into electronic plasma only to the
depth of $\lambda = c/\omega$ shows strong diamagnetic properties
of plasma. Here $c$ is the light velocity; $\omega$ is the Langmure
frequency (the frequency of "electronic sound");
$\omega = \sqrt {4\pi e^2n/m}$, where $n$ is the electronic density,
$m$ is the mass of electron, $e$ is its charge
\cite{Kad}. The value $\lambda = c/\omega$ plays the part of Debye
wavelength. In the paper of Arsenjev \cite{Ars} there is another indirect
evidence in favour of superconductivity. He indicates that collisions
in plasma leads to very small changes of momentum.

Let us try to describe a qualitative picture of the phenomenon under
the conjecture of superconductivity of plasma. Let us imagine ionized
plasma that generates a positively charged ball and a spherical
cover around this ball generated by negatively charged electronic
plasma. Suppose at last that in this cover it is realized
superconductivity. Electrons with momentum which tangent to the sphere
are deflected by the positive charge of the central part of the ball
remaining in the spherical cover. The superconductivity means that
the energy of electrons do not looses by collisions. Any continuous
vector field on a sphere has singularities. To avoid it, suppose that
in our case something like a geodesic flow on the sphere takes place.
Preservation of the ionization state means that one deals with the
cold plasma. This plasma does not failed due to superconductivity. It
explaines such a long life time of fireballs. A ball of plasma is
hanged in the air and in this ball takes place the current. It is
unusual current since it does not have a fixed direction at any point
and at any moment. It is spreaded in all possible tangent directions.
Inside the ball there is positive charge which directs and maintains
this current. Such phenomenon is possible at the state of
superconductivity only --- a current that does not have a direction.
Electronic liquid flows indeed, but it flows simultaneouly in all
directions. Hence the fireball is a physical realization of the mental
idea of Vlasov (which was used by him as the basis for his equation)
\cite{Vla}, the idea of probability distribution of electrons relative
to its positions and velocities.

How to organize such current at the low temperature under the usual
superconductivity? It would be interesting to take a negatively
charged supercoductor of the spherical form and to put inside it a
sufficiently large positive charge.

Thus the fireball, having the same density as the air and being of no
different from it except for its state, can freely float in the air
which is observed in reality. Its movement can obey airflows or
electromagnetic forces --- it is of no importance. What means the
explosion of fireballs? Some leakage of energy does exist, especially,
by contacts with charged bodies. Balance of energy which is necessary
to retain the current is failed, and plasma is extinguished. More
exactly, the charge of the kernel becomes unsufficient to maintain
electrons, they begin to fall into the kernel and this falling
decreases still further its charge. The process becames avalanche-type
and the fireball is quenched instantly. The pressure has no time to
level off. The difference in pressure stimulate a shock wave like that
of the spark or of the usual lightning. This wave is perceived as an
explosion. The known phenomenon --- the breakdown of the
superconductivity by sufficiently strong magnetic field --- explains
nonstandard interactions of fireballs with electrical equipment
\cite{Sta}.

The fireball is something like the model of the atom: a kernel and
electrons spreaded on orbits. Of course, we use the term "geodesic
flow" as a figurative expression. In accordance with the quantum
theory and with the Heisenberg's indeterminacy principle, electrons
are defined not as moving points but as a probability distribution.
Nevertheless the term "geodesic flow" reflects much more adequately
the bound state of electrons under the influence of powerful
electrostatic field of the positively charged kernel.

At latest works on the creation of fireballs \cite{Ego} it was
generated plasmoids (that was regarded by the authors as fireballs).
It has been found experimentally that: "It was formed the negatively
charged layer around the polarized kernel".

It may appear that the geodesic flow must lead to perpetual collisions of
electrons. But the volume of the cover is macroscopic and it is huge
by comparison with the size of an electron. Though any two big circles
of $S^2$ intersect but electrons moving along these orbits (no matter
how many they will be) can avoid collisions. Moreover, in the
approximation when electrons are considered as points, the probability of
collisions equals to zero. But the main reason: each separate
electron cannot scatter by himself being bounded coherently with other
electrons generating joint "plasma-sounded" vibrations. Thus the
approximation of the geodesic flow is reasonable.

Let us remark that all attempts to confine plasma was based on
ring-shaped regions. But it needs strong magnetic fields. May be it
would be better to try electrostatic fields. If there exists
superconductivity of plasma, then it must help. Attempts of confinment
are very difficult to realise because one endeavour to feed plasma by
energy, to clamp it down, to make a pinch. However we may give freedom
to plasma, to set it in a situation of "singing", when plasma will have
its own nondamping "sound" vibrations, moreover, vibrations which
resist attempts of attenuation. Then it would be a confinement and a
source of energy. It is exactly the state of superconductivity of
fireballs.

The Kapitza's conjecture on the nature of fireball was not exactly
correct, because his experiments with SHF-vibrations gave a
ring-shaped plasma rather than a ball. This was quite natural because
of the lack of the kernel. But the idea of feeding is very important.
One needs only firstly to create "atom", and after that to feed it by
energy.

One may imagine a mechanism of the fireball birth in the following way:
suppose, for instant, that the curl of positively charged particles
falls on the course of usual lightning. Then the negatively charged
plasma may surround this curl, and the fireball being stable begins
its life. This idea leaves room for experiments. For instance, to blow
a ball of ionized, positively charged gas with the electron beam, or
to lay it into the negatively charged medium. If it has been possible
to create fireballs artificially then an incorporation in it exact doses
of "fuel" and an energy-rejection would give a safe atomic reactor.

So far as in the picture of superconductivity of
Bardeen-Cooper-Schrieffer (BCS) the key part play phonons
(quasi-particles related to the fundamental frequences of the cristal
lattice), our goal is to find eigenfunctions of operator that
describes oscillations of spherical cover of plasma into the centrally
symmetrical electrostatic field and to prove the existence of a gap in
the energy spectrum which provides superconductivity.

Thus we take as a base the geodesic flow on the sphere $S^2$, that is
the uniform distribution of electrons on the bundle $T_*^1(S^2)$ of
tangent circles of radius 1 to the two-dimensional sphere of radius
$R$. The geodesic flow of electrons must be stable relative to small
disturbances, since, due to Coulomb interaction, electrons resist
compressions. Deviations from this stable movement provide
oscillations which play the role of phonons. Interaction of phonons
with electrons is the source of
superconductivity. One needs to interpret elementary vibrations around
the ground-state as independent oscillators and to replace it by
operators. As in an atom, the potential of the kernel exceeds the
interaction of electrons and it may be thought that electrons move
along great circles with the constant speed. That was the reason that
we begin with the bundle $T_*^1(S^2)$ of tangent circles to sphere.

The kernel of the fireball consists of positively charged ions. The
external cover of the fireball is the negatively charged plasma in
which there is a current --- the geodesic flow of electrons. That
means that the density of electrons is uniformly distributed on the
sphere; their velocities are uniformly distributed on the tangent
directions to the sphere and its absolute value equal to the speed of
circular motion around the kernel in response to its Coulomb
attraction. As a result, the potential of the external cover acting
upon internal points of the kernel equals to zero. This provide the
stability of the kernel. Indeed, if an ion penetrates into the cover
from the kernel then the electrostatic potential of the cover draw it
back, since the potential of the external (relative to this ion) part
of the cover equals to zero, meanwhile the potential of the internal
part of the cover draws it to the center. Electrons located outside of
the cover can not penetrate to the kernel, since the cover obstructs
this penetration. Namely, though the total charge of the fireball is
neutral, but in the vicinity of it, the negative charge of the cover
is prevailing that repuls closely-spaced negative charges preventing
its penetration into the kernel. Thus the positive charge of the
kernel remains practically unchanged. So, to explore the kernel one
can use the theory of nuclear matter \cite{Be}. Electrostatic
potential sustains the movement of electrons in the cover which
vaguely resembles the skin-effect in superconductors.

The fireball is something like soap-bubble. Soap-bubbles usually
slowly fall down in the air, because the pressure inside it (caused by
surface tension), and hence the density inside it, is a little more
than outside. So the weight of the bubble is a little greater than the
weight of the air. Similarly behave the fireballs. The additional
commplication is that the kernel tends to spred due to Coulomb forces,
and this electrostatic pressure is added to the usual pressure in gas.
The resulting pressure turns out to be a little greater than
atmospheric pressure which bring down the density in the central part
of the kernel. Meanwhile in the boundary part of the kernel (at the
vicinity of the cover) as well as in the cover itself the density
must be greater due to the pressure in the central part. As a result
the mean density of the fireball is only slightly greater than the
density of the air. Consequently the fireball can freely float in the air
with the slight tendency toward falling. The difference in pressure is
responsible for the effect of surface tension as well as for the
explosion of fireballs and the arising of shock waves when fireballs
are extinguished.

The main question remains to be answered. The question about losses of
the energy by scattering. Take as a pattern the BCS-theory of
superconductivity. In plasma of the cover there exist oscillations
like that of Bloch-waves in cristallic lattice. The difference is that
Bloch-waves are linked with the periodic potential of the lattice, while
collective oscillations of the cover are caused by the spherical
form of configuration space. The attracting field of the kernel and
the difference in the pressure inside and outside of the cover make
the cover resemble a spherical elastic film in which arise
spherical oscillations. Further on, the mechanism of plasma stability
becomes analogous to that of usual superconductivity: scattering of
electrons excite oscillations of the cover and that in its turn
stimulate electrons. This process can be regarded (in the quantum
level) as an exchanging interaction of electrons by using the analog
of phonons. Our aim is to find the energetic spectrum of the secondary
quantized problem and to detect a gap in this spectrum.

Let us recall that the notion of current in our case is different from
the usual electrodynamics notion. The movement of electrons is
uniformly distributed on all tangent directions. So in place of a
current one has only "microcurrents" which superimpose and do not give
electromagnetic fields. Consequently, the energy does not release: the
fireball does not radiate electromagnetic waves and there are no
losses of energy for that process. The fireball only fluoresce as usual
plasma.

\bigskip

     {\bf 1. Electronic cover}

\bigskip

Consider at first vibrations of the density of the geodesic flow on a
sphere. Imagine the cover as a homogeneous flux of a compressible fluid
on the bundle $T_*^1(S^2)$ of tangent circles at the two-dimensional
sphere of the radius $R$. On the sphere arise circular waves that can
be parametrized by its centers --- points of emittence of waves. At an
initial moment, take the disturbance of the density being the
$\delta$-function at a point $O$ with velocities uniformly
distributed on a tangent circle. As far as the velocity of superfluous
density is always directed outside, the wave bend the sphere passing
the equator, focusses to the diametrically opposite point $(-0)$
turning again into the $\delta$-function. We will call such waves
elementary. If a wave is not elementary then in place of focal points
appear another caustics. Each elementary wave on the bundle
$T_*^1(S^2)$ of tangent circles to the sphere is the torus $\T^2_0$.
Indeed in the initial moment one has the circle $S^1$ at the fiber
over the point $O$. Further, it is projected on a circle of the sphere
$S^2$ with the unically defined velocities which is orthogonal to the
circle. After two focusing, the initial circle is identically maps to
itself. Any two such torus intersect by a pair of anti-directed
geodesics, that connect centers of the corresponded waves.

Besides these waves of condensation there exist similar waves of
rarefaction. To quantize means to correlate an operator with each
elementary wave. Let us find the Hamiltonian for corresponding periodic
oscillations.

Note that the potential energy of the flow (let us denote it by $U$)
does not change during these oscillations. Indeed, from the point of
view of geodesic flow, the disturbances at each moment are
condenced at points of a circle which moves without deformation with
the constant speed along the torus $\T^2$. The degenerations of the
projection to a point in poles does not connect with the geodesic flow as
itself, but with its projection on the base. More exactly,
$T_*^1(S^2)$ --- the bundle of unit tangent vectors to $S^2$ --- is
diffeomorphic to the real three-dimensional projective space $\R\P^3$.
Tangent vectors to the geodesic flow are, say, left-invariant vector
field on $\R\P^3$ or on its two-fold covering $S^3$, that may be
realized as the group of unit in absolute value quaternions. Under the
stereographic projection from the center of the sphere $S^3$ on $\R^3$,
trajectories of this field transfer into a family of Clifford parallel
lines in $Ell(3)$ which is the space $\R^3$ equiped by elliptic metric
induced from $S^3$ by this projection. Shifts along the trajectories
of this family are isometries. To any fiber $S^1 \in T_*^1(S^2)$,
i.e. to the unit circle in
the tangent plane at a point $l$ of the base
$S^2$, correspond quaternions that commute with $l$: the circle in
$S^3$ passing through $l$ and the neutral quaternion. It is a greate
circle and so it projects onto a line in $Ell(3)$ passing through the
image of the neutral element. Under the Clifford shift, these lines
run along a toric surface. Hence, the set of trajectories of points
of a fiber over some point $x_0 \in S^2$ is two-dimensional torus. As
was shown by Clifford, the metric of this torus induced by the ambient
elliptic metric is Euclidian. Thus the potential energy of an
elementary wave remains constant.

Meanwhile, in contrast with the movement in a flat space, the momentum of
the wave changes. This is caused by the fact that the movement is
not free: it takes place not in the whole space but along the surface of
a sphere. Take the spherical coordinate on the sphere with the pole
$O$. Fronts of wave centered in $O$ are circles $\psi = Const.$ As an
momentum of this wave it should be taken the operator of
differentiation by $\psi$. But since this operator must be defined at
all points (in particular at poles), it is necessary to multiply it by
a function which is equal to zero at poles. The simplest such function
is $\cos \psi.$ Thus as an momentum we consider the operator

$$
p = i \frac {\hbar}{2}\cos \psi \frac {\partial }{\partial \psi}.
$$

We will use the rational system of units in which the light velocity
$c$ and the Plank's constant $\hbar$ equal to 1. So, the Hamiltonian is
defined by the operator

\beq
Hf =  - \cos \psi \frac {\partial }{\partial \psi}
(\cos \psi \frac {\partial }{\partial \psi})f + Uf,
\label{1}
\eeq
where $U$ up to multiplicative constant corresponds to potential
energy. We are interesting in solutions depending from $\psi$ only, so
$f = f(\psi)$. The equation for oscillations is

$$
- \cos^2 \psi \frac {d^2 f}{d \psi^2} +
\cos \psi \sin \psi \frac {d f}{d \psi} + Uf (\psi) = 0,
$$
or
\beq
L = \frac {d^2 f}{d \psi^2} -
\tan \psi \frac {d f}{d \psi} - \frac {U}{\cos^2 \psi}f = 0.
\label{2}
\eeq

We seek for the spectrum of the operator $L$ defined by the equation
(\ref{2}). Eigenvalues of $L$ will be denoted by $\alpha$. We obtain
the equation

\beq
(\cos^2 \psi)f'' - (\sin \psi \cos \psi)f'
+ (\alpha \cos^2 \psi - U)f = 0.
\label{3}
\eeq
Boundary conditions for the equation (\ref{3}) are

\beq
f \left ( \frac {\pi}{2} \right ) = f
\left ( -\frac {\pi}{2} \right ) = 0,
\label{4}
\eeq
since the function $f$ has to be independent on $\phi$ at poles.
Let us note that the boundary value problem
(\ref{3})-(\ref{4}) defines a selfajoint operator, because the equation
(\ref{3}) can be written in the form
$$
L_1(f)=\frac {d}{d\psi} \left ( (\cos \psi) \frac {df}{d\psi} \right )
- \frac {U}{\cos \psi}f = - \alpha \cos \psi f.
$$

Eigenfunctions of the selfajoint operator $L_1$ with the boundary
conditions (\ref{4}) relative to the positively defined selfajoint
operator of multiplication by the function $(\cos \psi)$ are mutually
orthogonal on the interval $[-\pi /2, \pi/2]$ with the weight
$(\cos \psi)$.

\bigskip

We found oscillations of density or the "longitudinal" oscillations on
the sphere. Let us find its "cross" oscillations, that is the
oscillations of the sphere itself from purely phenomenological point
of view. With this in mind, let us solve the following problem that has
an independent interest. It may be interpreted as a problem of
oscillation of the balloon cover.

Consider an elastic spherical film surrounding a volume with the
internal pressure greater than that of the outer air. We seek for its
eigen oscillations due tu surface tension.

On any element of the surface act forces of surface tension which are
applied to the boundary points of the element. The resultant of two forces
which act on the opposite points of a given direction at a point $O$
is normal to the surface and equals to the normal curvature of this
direction. The resultant of all forces is the normal to the surface
which equals to the integral of normal curvatures, hence, equals to the
mean curvature at a point $O$. For the stationary state, these forces are
balanced out by the difference in pressure. The
difference between rersultant and the force caused by the difference
in pressure (acting in the same direction) give the force acting on
the cover. As far as we know, the problem in this statement was not
considered. This model is different from usual statements of problems
in elastic theory since usually one considers forces to be
proportional to displacements.

Let us formalize the problem obtained. Introduce spherical coordinates
on the sphere of radius $R$.

$$
\{R\cos \psi \cos \phi , R\cos \psi \sin \phi , R\sin \psi \}.
$$

Denote the corresponding point of the unit sphere by

$$
r = \{\cos \psi \cos \phi , \cos \psi \sin \phi , \sin \psi \}.
$$

Denote a displacement of the point $R, \psi , \phi$ by

$$
R\Pi =
\{R\xi (\psi , \phi ), R\eta (\psi ,\phi ) , R\zeta (\psi \phi \}.
$$

The coordinates of the displaced point are $X = Rr + R\Pi$. Let us
introduce the moving frame on the unit sphere

$$
\begin{array}{l}
r = \{\cos \psi \cos \phi , \cos \psi \sin \phi , \sin \psi \}, \\
k = \{-\sin \phi , \cos \phi , 0\} = \frac {r_{\phi}}{\cos \psi}, \\
l = \{-\sin \psi \cos \phi , -\sin \psi \sin \phi , \cos \psi \}
= r_{\psi}.
\end{array}
$$

One has

$$
\begin{array}{l}
r_{\phi \phi} = (\cos \psi )k_{\phi}, \; k_{\phi} =
(-\cos \psi )r + (\sin \psi )l, \; k_{\psi} = 0, \; l_{\psi} = -r; \\
l_{\phi } = (-\sin \psi) k, \; r_{\phi \psi} = l_{\phi} =
(\sin \psi )k, \; r_{\psi \psi} = l_{\psi} = -r.
\end{array}
$$

Its vector products are

$$
[r \times k] = l, \; [k \times l] = r, \; [l \times r] = k.
$$

Calculate the mean curvature of the displaced surface, neglecting
members of second order relative to $\Pi$ and its derivatives

$$
X_{\phi} = (R\cos \psi)k + R\Pi_{\phi}, \;
X_{\psi} = Rl + R\Pi_{\psi}.
$$

The coefficients of the first quadratic form $ds^2$ equal

\beq
\begin{array}{l}
(X_{\phi},X_{\phi}) = R^2(\cos^2 \psi + 2\cos \psi (k,\Pi_{\phi})),\\
(X_{\phi},X_{\psi}) = R^2(\cos \psi (k,\Pi_{\psi}) + (l,\Pi_{\phi})),\\
(X_{\psi},X_{\psi}) = R^2(1 + 2(l,\Pi_{psi})).
\end{array}
\label{5}
\eeq

$$
[X_{\phi} \times X_{\psi}] = R^2 \{(\cos \psi )r +
(\cos \psi )[k \times \Pi_{\psi}] - [l \times \Pi_{\phi}]\}.
$$

$$
|[X_{\phi} \times X_{\psi}]| =
R^2 \sqrt {\cos^2 \psi + 2\cos^2 \psi (l,\psi) +
2\cos \psi (k,\Pi_{\phi})}.
$$

$$
\frac {1}{|[X_{\phi} \times X_{\psi}]|} =
\frac {1}{R^2\cos \psi} \{1 - (l,\Pi_{\psi}) -
\frac {(k,\Pi_{\phi})}{\cos \psi}\}.
$$

Normal to the surface is defined by the expression

\beq
n = r + [k \times \Pi_{\psi}] -
\frac {[l \times \Pi_{\phi}]}{\cos \psi} - (l,\Pi_{\psi})r -
\frac {(k,\Pi_{\phi})}{\cos \psi}r.
\label{6}
\eeq
We have

$$
\begin{array}{l}
X_{\phi \phi} = R\{(-\cos^2 \psi)r + (\cos \psi \sin \psi)l
+ \Pi_{\phi \phi}\}, \\
X_{\phi \psi} = R\{(-\sin \psi)k + \Pi_{\phi \psi}\},\\
X_{\psi \psi} = R\{-r + \Pi_{\psi \psi}\}.
\end{array}
$$
The coefficients of the second quadratic form
$dN^2$ equal
\beq
\begin{array}{l}
(X_{\phi \phi},n) = R\{-\cos^2 \psi - \cos \psi \sin \psi (r,\Pi_{\psi}))
+ (r,\Pi_{\phi \phi})\},\\
(X_{\phi \psi},n) = R\{\tan \psi (r,\Pi_{\phi})
+ (r,\Pi_{\phi \psi})\},   \\
(X_{\psi \psi},n) = R\{-1 + (r,\Pi_{\psi \psi})\}.
\end{array}
\label{7}
\eeq

Denote $(2 \times 2)$-matrix of the quadratic form $dN^2 - \lambda ds^2$
by $A(\lambda)$. Using (\ref{5}) and (\ref{7}) find entries of this
matrix

$$
\begin{array}{l}
A_{11} = R\{-\cos^2 \psi - \cos \psi \sin \psi (r,\Pi_{\psi})
+ (r,\Pi_{\phi \phi})\} - \lambda R^2 \{\cos^2 \psi +
2\cos \psi (k,\Pi_{\phi})\},\\
A_{12} = R\{\tan \psi (r,\Pi_{\phi}) + (r,\Pi_{\phi \psi})\} -
\lambda R^2\{\cos \psi (k,\Pi_{\psi}) + (l,\Pi_{\phi})\},\\
A_{22} = R\{-1 +(r,\Pi_{\psi \psi})\} -
\lambda R^2\{1 + 2(l,\Pi_{\psi})\}.
\end{array}
$$

Roots of the equation $det A = 0$ define main curvatures. By dividing
half of the coefficient of $\lambda$ in this equation by the first
coefficient and neglecting members of the second order relative to
$\Pi$ we find the mean curvature

$$
K = \frac {1}{R}
\left \{
-1 - \frac {\tan \psi}{2} (r,\Pi_{\psi})
+ \frac {(r,\Pi_{\phi \phi})}{2\cos^2 \psi} +
\frac {(r,\Pi_{\psi \psi})}{2} + (l,\Pi_{\psi}) +
\frac {(k,\Pi_{\phi})}{\cos \psi}
\right \}.
$$

The force which act on the surface at a point $(\phi, \psi)$
is directed along the normal $n$ (see (\ref{6})). It is proportional
to $K$ with the coefficient of propotionality which we denote by
$\beta$. The inner pressure gives the force which is proportional to
$\frac {1}{R}n$. To have an equilibrium for the case when the
displacement of surface equls to zero (sphere), the coefficient of
proportionality must be also equal to $\beta$. The summary force is
equal

$$
\beta r
\left \{
- \frac {\tan \psi}{2}(r,\Pi_{\psi}) +
\frac {(r,\Pi_{\phi \phi})}{2\cos^2 \psi} + \frac {r,\Pi_{\psi \psi}}{2}
+ (l,\Pi_{\psi}) + \frac {(k,\Pi_{\phi})}{\cos \psi}
\right \}.
$$

One obtains the equation

\beq
R \frac {\partial^2 \Pi}{\partial t^2} =
\beta r
\left \{
\{-\frac {\tan \psi}{2}(r,\Pi_{\psi}) +
\frac {(r,\Pi_{\phi \phi})}{2\cos^2 \psi} +
\frac {(r,\Pi_{\psi \psi})}{2} + (l,\Pi_{\psi}) +
\frac {(k,\Pi_{\phi})}{\cos \psi}
\right \}.
\label{8}
\eeq

Rewrite the equation (\ref{8}) in the moving frame coordinates:
$u=(r,\Pi), \, v=(k,\Pi), \, w=(l,\Pi).$

$$
\begin{array}{l}
\frac {1} {\beta} \frac {\partial^2 u}{\partial t^2} =
\frac {1}{R}
\left \{
- \frac {\tan \psi}{2}u_{\psi} +
\frac {u_{\phi \phi} + (\cos^2 \psi) u_{\psi \psi}}{2\cos^2 \psi} +
w_{\psi} + \frac {v_{\phi}}{\cos \psi}
\right \}. \\
\frac {1} {\beta} \frac {\partial^2 v}{\partial t^2} =
(k,\frac {\partial^2 \Pi}{\partial t^2}) = 0, \\
\frac {\partial^2 w}{\partial t^2} =
(l,\frac {\partial^2 \Pi}{\partial t^2}) = 0.
\end{array}
$$

If at the initial moment
$v=0, \, \dot v = 0$ and $w=0, \, \dot w = 0$, then
$v \equiv 0, \, w \equiv 0$. Thus, the linearization leads to purely
radial vibrations. The final equation is

\beq
\frac {1} {\beta} \frac {\partial^2 u}{\partial t^2} =
\frac {1}{R}
\left \{
-\frac {\tan \psi}{2}u_{\psi} +
\frac {u_{\phi \phi} + (\cos^2 \psi) u_{\psi \psi}}{2\cos^2 \psi}
\right \}.
\label{9}
\eeq

\smallskip

As expected, the right hand side of the equation (\ref{9}) gives the
Laplace-Beltrami operator on the sphere.

As soon as the equation (\ref{9}) is autonomous relative to
$\phi$, we find the solution in the form of traveling wave relative to
$\phi$, that is $u = u(\phi - \kappa t, \psi)$. We have

$$
\frac {\kappa^2}{\beta} u_{\phi \phi} = \frac {1}{R}
\left \{
- \frac {\tan \psi}{2}u_{\psi}
+ \frac {u_{\phi \phi}}{2\cos^2 \psi} +
\frac {u_{\psi \psi}}{2}
\right \}.
$$

Separation of variables $u = \Phi (\phi) \Psi (\psi)$ gives

$$
\frac {\Phi''}{\Phi} =
\frac{(-\sin \psi \cos \psi) \frac {\Psi'}{\Psi} +
(\cos^2)\psi\frac {\Psi''}{\Psi}}{2R\kappa^2\cos^2 \psi - \beta} = C.
$$

The period of the function
$\Phi$ must be equal to $2\pi$, hence
$C = - n^2, \; n \in \Z$. The equation for $\Psi$ has the form

\beq
(\cos^2 \psi)\Psi'' - (\sin \psi \cos \psi)\Psi' +
(2n^2R\kappa^2 \cos^2 \psi - n^2\beta)\Psi = 0.
\label{10}
\eeq

The boundary conditions for the equation
(\ref{10}) are

\beq
\Psi \left ( \frac {\pi}{2} \right ) =
\Psi \left ( -\frac {\pi}{2} \right ) = 0,
\label{11}
\eeq
since the function $\Psi$ does not depend on $\phi$ at poles. The
expression $\alpha =2n^2R\kappa^2$ can be considered as the spectral
parameter for the problem (\ref{10})-(\ref{11}). Zonal spherical
harmonics are its solutions.

\smallskip

The spectral problem
(\ref{10})-(\ref{11}) appears exactly the same as
(\ref{3})-(\ref{4}). The striking fact that entirely different
approaches and different phenomena lead to the same equation
(\ref{10})=(\ref{3}) give confidence to the Hamiltonian
(\ref{1}).

\bigskip

Eigenfunctions $u_{m,l}$ of the Laplace-Beltrami operator on the
two-dimensional sphere are enumerated by pairs of integers
$(m,l),\; |m| \le l$ and relate to the eigenvalues $- l(l+1)$.

In accordance with the established linguistic tradition, we shall call
particles which are obtained by quantization of the above-mensioned
waves (both langitudinal and cross waves) by spherons.

\bigskip

\begin{theorem}
The solution $u=0$ of the equation (\ref{9}) is stable in Liapunov's
sense.
\end{theorem}

\Proof. The equation (\ref{9}) can be rewritten as the ordinary
differential equation in the space $l_2$

\beq
\frac {\partial^2 u}{\partial t^2} = Au.
\label{12}
\eeq

The spectrum of the operator $A$, which stand at the right hand side
of the equation (\ref{12}), equals

$$
\lambda_{m,l} = - l(l+1), \quad l \ge |m|.
$$

The equation (\ref{12}) can be written in normal form as a system

\beq
\begin{array}{l}
\frac {\partial u}{\partial t} = v, \\
\frac {\partial v}{\partial t} = Au.
\end{array}
\label{13}
\eeq
The spectrum of the matrix

$$
\left (
\begin{array}{cc}
0 & I \\
A & 0
\end{array}
\right ),
$$
standing at the right hand side of the system (\ref{13}),
is purely imaginary. It is equal

$$
\pm i\sqrt{l(l+1)}.
$$

Indeed, if $f_k, \; k \in \N$ is a basis of eigenvectors of the
operator $A$ responding to eigenvalues $-{\lambda_k^2}$, then

$$
u = f_k e^{\pm i\lambda_k t},
v = {\pm i\lambda_k}f_k e^{\pm i\lambda_k t}
$$
is a fundamental system of solutions to the equation (\ref{13}).

The theorem is proved.

$\Box$

\bigskip

{\bf 2. Second quantization}

\bigskip

We will quantize the equation obtained in momentum space. The role of the
momentum space will plays the space of spherical harmonics. Let us note
that it has the discrete basis that may be defined by eigenfunctions
of the Laplace-Beltrami operator. The thechnique of quasi-discrete
representation, i.e. the preliminary consideration of the problem
in finite parallelepiped with periodic boundary conditions and after
that the passage to the limit from discrete case to the continuous
momentum representation \cite{Bog} (one of the feeble item in the ground of
quantum statistics), is found to be unnecessary due to compactness of
the phase state.

Denote the spherical harmonics by $v_{m,l}(s),\; s \in S^2$.
Expand a solution to the equation

\beq
\frac {\partial^2 u}{\partial t^2} = \Delta u,
\label{14}
\eeq
where $\Delta$ is the Laplace-Beltrami operator on the sphere $S^2$,
in the series in terms of spherical harmonics
$v_{m,l}(s)$. Coefficients of this series are defined by the
expressions

\beq
\hat u_{m,l}(t) = \int_{S^2} u(s,t)v_{m,l}(s)\, ds.
\label{15}
\eeq

Apply the operator $\frac {\partial^2}{\partial t^2}$ to the
(\ref{15}) taking into account that the Laplace-Beltrami operator is
selfajoint and that the functions $v_{m,l}(s)$ are its eigenfunctions
with the eigenvalues $- l(l+1)$

$$
\frac {\partial^2 \hat u_{m,l}(t)}{\partial t^2} =
\int_{S^2} \Delta u(s,t) \, v_{m,l}(s) ds =
\int_{S^2} u(s,t) \, \Delta v_{m,l}(s) ds = - l(l+1)\hat u_{m,l}(t).
$$

One obtains a decomposed infinite system of ordinary differential
equations

\beq
\frac {\partial^2 \hat u_{m,l}(t)}{\partial t^2} = -
l(l+1)\hat u_{m,l}(t).
\label{16}
\eeq

Each equation of the system (\ref{16}) is the equation of harmonic
oscillator with the unit mass and the frequence $\sqrt{l(l+1)}$.
We shall call it by the oscillator $(m,l)$. Apply to it the standard
proceedure of quantization. That leads to second quantization of the
oscillation of the geodesic flow of electrons.

The Hamiltonian of the oscillator $(m,l)$ is

$$
H^0_{m,l} = \frac {1}{2} (p^2 + l(l+1)q^2).
$$

We enumerate eigenfunctions of the operator $H^0_{m,l}$ by indices
$n$ and denote it by $\xi_n (m,l)$ (the number of an eigenfunction is
denoted by the lower index). It is well known that the corresponding
eigenvalues are
$E_n = (n + \frac {1}{2})\sqrt{l(l+1)}$. Introduce the
operators \cite{Wei} acting on eigenfunction as follows

\beq
a^+(m,l) \xi_n(m,l)  = \sqrt{n + 1} \, \xi_{n+1}(m,l), \quad
a^-(m,l) \xi_n(m,l)  = \sqrt{n} \, \xi_{n-1}(m,l).
\label{17}
\eeq

The ground-state of the system corresponds to $n = 0$, since in this
case the next conditions are fullfiled $a^-(m,l)\xi_0(m,l)=0$, and

\beq
\xi_n(m,l)  = \frac {(a^+(m,l))^n}{\sqrt{n !}} \xi_0(m,l).
\label{18}
\eeq

The formulas (\ref{17}) show that states $\xi_n(m,l)$ are
eigenfunctions of the operators $a^-(m,l) a^+(m,l)$ and
$a^+(m,l) a^-(m,l)$ with the eigenvalues $n + 1$ and $n$
correspondingly. Consequently, it is natural to call the operator
$a^+(m,l) a^-(m,l)$ by the operator of number of spherons
of type $(m,l)$. Denote this operator by $\hat n^s (m,l).$ The
commutation relations take the standard form

$$
[a^-(m,l),a^+(m,l)] = 1, \; [a^-(m,l),a^-(m,l)] =
[a^+(m,l),a^+(m,l)] = 0.
$$

Meanwhile operators that relate to noncoinsiding values of indices
commute. The Hamiltonian for an oscillator is $H^0(m,l) =
\sqrt{l(l+1)}(a^+(m,l)a^-(m,l) + \frac {1}{2}{\rm I}^s(m,l))$, or

$$
H^0(m,l) = \sqrt{l(l+1)}(\hat n^s (m,l) + \frac {1}{2}{\rm I}^s(m,l)),
$$
where ${\rm I}^s(m,l)$ is the identity operator on the space of states
of the oscillator $(m,l)$. The full Hamiltonian of the field of the
cross waves is represented by the sum

\beq
H^0 = \sum_{l,|m| \le l} H^0(m,l) =
\sum_{l,|m| \le l} \sqrt{l(l+1)}(\hat n^s (m,l) +
\frac {1}{2}{\rm I}^s(m,l)).
\label{19}
\eeq

One can express the Hamiltonian in the representation of occupation
numbers. For brevity sake, let us relabel Fourier coefficients
$(m,l)$ of spherical harmonics by unique index $h$ and denote the
corresponding operators by $a^+(h)$ and $a^-(h)$. The state of the
system is characterized by a set of occupation numbers
$h_1, h_2, ... h_N$, i.e. the state

$$
\xi (h_1, ... h_N) = \prod_{1 \le k \le N}
\frac {(a^+(k))^{h_k}}{\sqrt{h_k !}} \xi_0(k)
$$
means that there are $h_i$ particles of the sort $\xi_i$.

The mechanism of the superconductivity of plasma is somewhat different
from that of electron-phonon interaction in BCS-theory. The later is
connected with interactions of electrons and vibrations of ion lattice,
that is with something external relative to the electron flow itself.
Whereas electron-spheron interaction is the interaction of electrons
with vibrations of the flow of electrons. In other words, this is
something like the selfconsistency of Hartree \cite{Har}. The only
difference is that in Hartree theory the movement of an electron is
consistent with its own field, wheareas in plasma, the movement of an
electron which affect the creation of the field of the flow must be
consistent with the field of the flow as a whole. In this situation
the field of the flow is not the electromagnetic field, but the field
of spheron's oscillations. Roughly speeking, fluctuations in state of
an electron as if "damp" the electron flow giving to it quants of
energy. In respons, the coherent flow slightly "pushes" this electron
returning it back into "formation".

As the components of the field of the flow it may be taken longitudinal and
cross waves. The equation (\ref{3}) describes longitudinal
oscillations in the absence of cross ones; the equation
(\ref{10}) describes cross oscillations in the absence of longitudinal
ones. But in the general case longitudinal and cross oscillations
are mutually connected. Local deformations of the sphere lead to the
changes of density; local decreasings of density decrease local
forces of surface tension that leads to deformation of the surface.
Note that the increasing of the function $f$ decrease the function
$\Psi$ and vice versa. In other words oscillations of the longitudinal
and cross waves are in opposite phase. Taking into account that these
oscillations satisfy to the same equation, it will be supposed that
$\Psi = - f$. So spherons will be described by a scalar field with the
Hamiltonian (\ref{19}). Note further, that an elementary wave of
density arrive at the diametrically opposite point being in the same
phase as at the initial point, that is the excess of density at the
initial point leads to the excess of density after focussing. Thus we
will be interested only in even functions relative to the antipodal
mapping. Hence, the spectral parameter $l$ must be an even number. The
corresponding system of eigenfunctions are complete in the projective
space $\R\P^2$. The Hamiltonian of the field of spherons takes the form

\beq
H^0 = \sum_{l,|m| \le 2l} \sqrt{2l(2l+1)}(\hat n^s (m,2l) +
\frac {1}{2}{\rm I}^s(m,2l)).
\label{20}
\eeq

Our aim is to write the Hamiltonian which consists of

1). the free Hamiltonian of the spherons field (\ref{20}),

2). the free Hamiltonian of the field of electrons, and

3). the interaction Hamiltonian.

\bigskip

{\bf 3. The gap}

\bigskip

We wish to reduce the situation to the case of the Bose statistics as
in usual theory of superconductivity. In that theory it was used
quasi-particles --- Cooper's pairs. There is some mystery in the fact
that a pair consists of electrons with the opposite momentum and
spins. It is rather strange that particles flying in opposite
directions are linked. Namely Cooper that proposes this idea calls this
connection "deeply mysterious". Though, there was some statistical
arguments in favour of it (see, for instant, \cite{Schm}). The main
reason is that the combined spin becomes zero which leads to the
Bose statistics. So in the state with minimum energy there are plenty
of particles, so called coherent condensate, having high-capacity
inertia which obstacle small changes. Together with the energetic gap
this gives the intuitive explanation of the superconductivity.

In our case the combination in Cooper's pairs is much more natural
than in BCS-theory. We again combine electrons with the opposite spin
and momentum. But now electrons of the pair move along the same
geodesic in the same direction. It is convenient to think on it as
situated at the diametrically opposite points. More exactly, the support
of the function that describes the joint state of the pair belongs
to the antipodal diagonal of the direct product $S^2 \times S^2$, that
is it belongs to the manifold $\{x,-x\} \in S^2 \times S^2$. In
particular this means that if one of the electrons was detected at the
point $x$ then its pair will be at $-x$. Thus if one of the electrons
of the pair jumps from the state with the momentum $p$ into that with the
momentum $k$ then the wave of density having the concentration at
diametrically opposite points transfers its pair
from the state with the momentum $-p$ into that with the momentum $-k$.

Let us show that combining electrons in Cooper's pair provides the gap
in energy spectrum. Denote the width of the electronic cover by $r_0.$
The phase space for electronic cover related to spheres with radii
from $R$ to $(R+r_0)$ is the direct product of the manifold of
geodesics on a sphere (that is $\R\P^2$) by the interval $[R, R+r_0]$.
An momentum $k$ and a spin, taking values $\pm 1/2$, parametry
eigenfunctions. The set of eigenfunctions plays the part of a
crystalline lattice being something like a framework maintaining
superconductivity. Accidental disturbance in density of electrons
spreads along the sphere as a spheron and compensates the related
disturbance. It is possible since the phase velocity of waves is
far much than that of electrons. The kinetic energy of an electron is
defined by the radius of the sphere on which it moves. Aside from
interchanging processes due to spherons, the potential energy is
defined by interelectron repulsion which we will often neglect in what
follows. Interchanging processes give the effective attraction, and
this attraction is far-acting since oscillations of electronic cover
is a global process acting on the whole sphere. So the interchanging
processes give a negative contribution into potential energy of the
system. That is the physical reasson for arising of the energetic gap in
the spectrum. Let us note that the role of Debye wavelength which
defines the characteristic linear length of the system, plays the width
of the electronic cover $r_0$.

Denote by $Q$ the charge of the kernel of a fireball, and by $R$ its
radius. Let $\rho$ be the average volume density of ions, and $e$ its
charge (in the case of single ionization). One has
$Q = \frac {4}{3}\pi R^3 \rho e$. Let $\mu$ be the average volume
density of electrons in the cover. The volume between spheres of the
radii $R$ and $R+r$ equals

$$
\frac {4}{3}\pi (R+r)^3 - \frac {4}{3}\pi R^3.
$$

The attractive force to the kernel acting on an electron that lie at
the distance $R+r$ from the center of the kernel equals
$\frac {Qe}{(R+r)^2} = \frac {4\pi \rho R^3 e^2}{3 (R+r)^2}$. But this
force is shielded by the negative charge of the layer between radii
$R$ and $R+r$. Denote the shielding coefficient (depended on $r$ only)
by $\kappa (r)$. Then the shielded attractive force of the kernel
equals $\kappa (r)\frac {4\pi \rho R^3 e^2}{3 (R+r)^2}$. Besides this
force, on an electron acts the force of repulsion from electrons of
the layer $(R,R+r)$ only, since the resultant force acting from the
external layer equals to zero. The resultant force acting from the
internal layer equals to that of repulsion from the center induced by
the overall charge of the layer $(R,R+r)$. Hence, the full force is
equal to

$$
F = \kappa (r)\frac {4\pi \rho R^3 e^2}{3(R+r)^2} -
\left (
\frac {4\pi (R+r)^3}{3(R+r)^2} - \frac {4\pi R^3}{3(R+r)^2}
\right )
\mu e^2.
$$

This force equals to zero at $r=r_0$, consequently, the equation defining
$r_0$ is

\beq
\kappa (r_0) \rho R^3 = ((R+r_0)^3 - R^3)\mu.
\label{21}
\eeq

The equation (\ref{21}) is in a sense the condition of neutrality of
the fireball. If one denote $\frac {r_0}{R} = \delta$, then the
equation for $\delta$ is

\beq
\delta^3 + 3\delta^2 + 3\delta = \frac {\kappa (\delta) \rho}{\mu}.
\label{22}
\eeq

Find the velocity $v$ of electron moving at the orbit of radius
$R+r$. The equation of movement at the plane $(x,y)$ are

$$
x = (R+r) \cos \frac {v}{R+r}t; \qquad  y = (R+r) \sin \frac {v}{R+r}t.
$$

By differentiation one obtains the centrifugal acceleration

$$
a = \frac {F}{m} = \frac {v^2}{R+r},
$$
where $m$ is the mass of the electron. Hence,

$$
v = \sqrt {\frac {F(R+r)}{m}}.
$$

The kinetic energy of an electron is

$$
\varepsilon (R+r) = \frac {mv^2}{2} =
\frac {2\pi e^2}{3(R+r)}[\kappa (r)\rho R^3 - \mu ((R+r)^3 - R^3)].
$$

It follows that

$$
\varepsilon (R) = \frac {2\pi \rho e^2}{3} R^2.
$$

Besides, in view of (\ref{22}), $\varepsilon (R+r_0) = 0.$

As in \cite{Sch},
imagine all admissible states in which electrons with momentum and spin
$((p,1/2),\; (-p,-1/2))$ are combined in pairs. The scattering of a paire
$((p,1/2),\; (-p,-1/2))$ on a pair $((k,1/2),\; (-k,-1/2))$ means the
transition from the state $\psi$ in which the cell
$((k,1/2),\; (-k,-1/2))$ is free and the cell $((p,1/2),\; (-p,-1/2))$
is occupied, into the state $\chi$ being in opposite condition. Denote
by $x_p^2$ the probability that a cell $((p,1/2),\; (-p,-1/2))$
is occupied. Then the probability that a cell $((k,1/2),\; (-k,-1/2))$
is free equals $y_k^2 = 1-x_k^2$. The probability of the state $\psi$ in
which the cell $((k,1/2),\; (-k,-1/2))$ is free and the cell
$((p,1/2),\; (-p,-1/2))$ is occupied is the product of probabilities
$P_{\psi}^2 = x_p^2(1-x_k^2)$. The amplitude of the state $\psi$ is
$P_{\psi} = x_p y_k$. Analoguously the amplitude of the state $\chi$
in which the cell $((p,1/2),\; (-p,-1/2))$ is free and the cell
$((k,1/2),\; (-k,-1/2))$ is occupied equals $P_{\chi} = x_k y_p$. Define
the energy of an electron with the momentum $p$ by $\varepsilon (p)$, and
denote the matrix element of the potential energy of interaction for
the process in question by $W_{p,k}$. The full energy of electronic
cover (discarding the Coulomb forces of electronic repulsion) is

\beq
E = \sum_p 2\varepsilon (p) x_p^2 +
    \sum_{p,k} W_{p,k} x_p x_k y_p y_k =
    \sum_p 2\varepsilon (p) x_p^2 -
    W(\sum_p  x_p y_p)(\sum_k  x_k y_k).
\label{23}
\eeq
Here $W_{p,k}$ is changed by $-W$ for momentum related to tangent
bundle to spheres of radii in the interval
$(R,R + r_0)$, and $W_{p,k} = 0$ for other momentum. The factor 2 of
the first sum arised from the fact that electrons in the states both
$(p,1/2)$ and $(p,-1/2)$ have the energy $\varepsilon (p)$.

The ground-state gives the minimum to the energy $E$. Let us find this
minimum.

\bigskip

\begin{theorem}
The following relation is fulfilled at the point of minimum of the
energy $E$

$$
1 = \sum_l \frac {2W}{\sqrt {\Delta^2 + \varepsilon^2(l)}},
$$
where
$$
\Delta = W \sum_k x_k y_k = W\sum_k \sqrt{x_k^2 - x_k^4} .
$$
\end{theorem}

\Proof.

Differentiate the expression for $E$ relative to $x_l^2$

$$
\frac {\partial E}{\partial (x_l^2)} =
2\varepsilon (l)  - W (\sum_k x_k y_k) \frac {1-2x_l^2}{2x_l y_l} = 0.
$$
The notation used in the statement of the theorem:
\beq
\Delta = W \sum_k x_k y_k = W\sum_k \sqrt{x_k^2 - x_k^4}
\label{24}
\eeq
gives the possibility to rewrite the preceding equation
in the form
\beq
2\varepsilon (l) \sqrt {x_l^2 - x_l^4} = \Delta (1 - 2x_l^2).
\label{25}
\eeq

By squaring we obtain the quadratic equation relative to $x_l^2$.

$$
x_l^4 - x_l^2 + \frac {\Delta^2}{4(\Delta^2 + \varepsilon^2 (l))} = 0.
$$

$$
x_l^2 = \frac {1}{2} \left (
1 \pm \sqrt {\frac {\varepsilon^2 (l)}{\Delta^2 +
\varepsilon^2 (l)}}\right ).
$$

The sign minus should be taken in this expression, since otherwise the
right hand side of the equation (\ref{25}) would be negative. So

$$
y_l^2 = \frac {1}{2}
\left (
1 + \sqrt {\frac {\varepsilon^2 (l)}{\Delta^2 + \varepsilon^2 (l)}}
\right ).
$$

We shall show below that, as in the classical theory of superconductivity,
the value $\Delta$ corresponds to the gap in the energy spectrum of
the system. Substituting this expression into the formula
(\ref{24}), we obtain the analog of the famous gap equation:
$$
\Delta =  \sum_l \frac {W}{2}
\left (
1 - \frac {\varepsilon^2(l)}{\Delta^2 + \varepsilon^2(l)}
\right )^{1/2}
$$
or
\beq
1 = \sum_l \frac {2W}{\sqrt {\Delta^2 + \varepsilon^2(l)}}.
\label{a}
\eeq

The theorem is proved.

$\Box$

\bigskip

Let us introduce $\nu (\varepsilon)$ --- the density of states of
electrons with the energy $\varepsilon$. Using this notion rewrite the
gap equation in an integral form

$$
1 = W  \int_{\varepsilon (R+r_0)}^{\varepsilon (R)}
\frac {\nu (\varepsilon)d\varepsilon}{2\sqrt {\Delta^2 + \varepsilon^2}}.
$$

The oscillation of the function $\nu (\varepsilon)$ on the interval
$(R, R+r_0)$ is insignificant. We shall denote $\nu (\varepsilon)$
simply by $\nu$ and it can be carried out of the sign of the integral.
We obtain

\beq
\frac {2}{W\nu} =   \int_{\varepsilon (R+r_0)}^{\varepsilon (R)}
\frac {d\varepsilon}{\sqrt {\Delta^2 + \varepsilon^2}}.
\label{26}
\eeq

The integral at the right hand side of (\ref{26}) monotonically
decreases from $+\infty$ to $0$ as $\Delta$ runs from $0$ to
$+\infty$. Consequently, a solution $\Delta$ to the equation
(\ref{26}) exists and is unique. Let us find the explicit expression
for $\Delta$ by calculating the integral (\ref{26}).

$$
\frac {2}{W\nu} =  \mathop{\rm Arsh}
\frac {\varepsilon (R)}{\Delta} -
\mathop{\rm Arsh}\frac {\varepsilon (R+r_0)}{\Delta} =
\mathop{\rm Arsh}\frac {2\pi e^2 \rho R^2}{3\Delta}.
$$

Hence,

\beq
\mathop{\rm sh} \frac {2}{W\nu} = \frac {2\pi e^2 \rho R^2}{3\Delta}.
\label{27}
\eeq

It follows

\beq
\Delta = \frac {2\pi e^2 \rho R^2}{3} \mathop{\rm sh^{-1}}
\frac {1}{2W\nu}.
\label{28}
\eeq

Let us show that $\Delta$ is actually the gap in the energy
spectrum. Having this in mind, we evaluate the contribution $\xi (l)$
from the adding of the combined pair of electrons with momentum
$\pm l$ into the energy of the ground-state (when all electrons are
combined in pairs). From (\ref{23}) we have

$$
\xi (l) = 2\varepsilon (l) x_l^2 -2Wx_ly_l\sum_k x_ky_k.
$$

Indeed, the first summond in this formula is the kinetic energy of the
pair $(l,-l)$, and the second is the contribution into the negative
part of the energy of the ground-state arising from the possibility of
interaction processes of the pair in question with all other pairs
$(k,-k)$. The coefficient 2 reflects the fact that the pair $(l,-l)$
in the sum (\ref{23}) occurs twice: in summing by $k$, and in summing
by $p$. Itroduce the notation
$E_l = \sqrt {\Delta^2 + \varepsilon^2 (l)}.$ Taking into account
values: $x_l^2 = \frac {1}{2} (1 - \frac {\varepsilon (l)}{E_l})$, and
denotation $\Delta = W\sum_k x_ky_k$, we find

$$
\begin{array}{c}
\xi (l) = 2\varepsilon (l)\frac {1}{2}(1 -
\frac {\varepsilon (l)}{E_l}) - 2\Delta
\left (
\frac {1}{4}(1 - \frac {\varepsilon^2 (l)}{E_l^2})
\right )^{1/2} = \\
\varepsilon (l) - \frac {\varepsilon^2 (l)}{E_l} -
\Delta \frac {\Delta}{E_l} = \varepsilon (l) - E_l.
\end{array}
$$

If to add to the ground-state merely one (unpaired) electron with the
momentum $l$ then the pairs $(l, -l)$ cannot contribute into the energy
of the ground-state, and hence the energy of the system after such
adding become

$$
E - \xi (l) + \varepsilon (l) = E + E_l = E +
\sqrt {\Delta^2 + \varepsilon^2 (l)}.
$$

Consequently, the adding of an isolated electron raises the energy of the
ground-state at least on $\Delta$. This means that to change the
ground-state, i.e. to disjoin a combined pair, it is required the
energy non less than $2\Delta$. This is the same as saying that there
is a gap in the energy spectrum.

\bigskip

To estimate the shielding effect $\kappa (r)$, we will use the
methodology of Tomas-Fermi, developed for computations of heavy metals
spectra. The Tomas-Fermi model for centrally symmetric fields is based
on the supposition that there are greate many electrons in a domain
where the potential changes slowly, and in this domain the Fermi
statistics is applicable. This supposition is fulfilled in our model
of fireballs.

\bigskip

Firstly, let us remined the deduction of the Tomas-Fermi equation in
suitable denotations. Denote by $p_0(r)$ the maximal absolute value of
momentum of electrons situated at the distance $r$ from the center.
Let us agree that $\frac {p_0^2(r)}{2m}$ is the excess of the energy
of an electron over the level of the termal potential. We have

\beq
V(r) + \frac {p_0^2(r)}{2m} = 0.
\label{29}
\eeq

Let an electron moves freely in a domain $\Omega$ with the electron
density $\mu = \frac {N}{\Omega}$. Then the number of quantum states
with the absolute values of momentum from $p$ to $p+dp$ is

$$
2\frac {\Omega}{(2\pi)^3}4\pi k^2 dk,
$$
where $\hbar k = p$, and the factor 2 takes into account two possible
orientations of the spin. We integrate from 0 to $k_0$ and equate the
result to the full number of electrons $N$:

$$
2\frac {\Omega}{(2\pi)^3} \frac {4\pi}{3}k_0^3 = N,
$$
or
$$
k_0^3 = 3\pi ^2 \mu.
$$

Let now suppose that a subdomain $O \subset \Omega$ is selected being
sufficiently large for correctness of (\ref{29}) and sufficiently
small, such that the potential energy in it changes not too much. Then

\beq
\mu = \frac {(2m)^{3/2}}{3\pi^2\hbar^3} (-V^{3/2}).
\label{30}
\eeq

By substitution (\ref{30}) into the Poisson equation for electrostatic
potential $-V/e$ with the density of charge $- e\mu$:
$$
\Delta V = -4\pi e^2 \mu,
$$
one obtaines
$$
\frac {1}{r} \frac {d^2}{dr^2} (rV(r)) = -
\frac {4e^2}{3\pi \hbar^3} (2m)^{3/2} (-V)^{3/2}.
$$

Using dimensionless variables
$$
x = \frac {r}{b} = \frac {2^{7/3}}{(3\pi)^{2/3}} r; \quad
f = - \frac {rV}{Qe^2},
$$
where $Q$ is the total charge of the kernel, one obtains the Tomas-Fermi
equation

\beq
\frac {d^2 f}{dx^2} = \frac {f^{3/2}}{\sqrt {x}}.
\label{31}
\eeq

\bigskip

Let us apply the equation (\ref{31}) to our theory of fireballs. While
calculating heavy metal spectra, the following boundary values for the
equation (\ref{31}) is customarily used (cf. \cite{Sch}, \cite{Bet})

\beq
f(0) = 1; \quad \lim_{x \to \infty} f(x) = 0.
\label{32}
\eeq

However, we need for our model another boundary conditions. It must be
taken that at the distant $R$ we have the Coulomb potential, and the
cover ends at points where termal fluctuations of
electrons majorize the attractive forces. At the boundary points, the
potential will be equal to the fixed value $- \zeta$, where $\zeta$
corresponds to the termal energy of electrons at the temperature $T$.
The related boundary conditions are

\beq
f(R) = 1; \quad
f(R+r_0) = \frac {R+r_0}{Qe^2} \zeta.
\label{33}
\eeq

Let us explore the geometrical structure of solutions to the equation
(\ref{31}).

\begin{theorem}
The equation (\ref{31}) admits the group $G$ of scale symmetries
allowing to reduce it to a two-dimensional system which does not
depend on the argument $x$.
\end{theorem}

\Proof.

Let us rearrange the equation (\ref{31}) as the system

\beq
\begin{array}{l}
\frac {df}{dx} = g \\
\frac {dg}{dx} = \frac {f^{3/2}}{\sqrt {x}}.
\end{array}
\label{34}
\eeq

The group $G$ can be found in the form
$f \mapsto \lambda f, \;
g \mapsto \mu g, \; x \mapsto \nu x$. These transformations
transfer the system (\ref{34}) into the system

$$
\begin{array}{l}
\frac {\lambda}{\nu}\frac {df}{dx} = \nu g \\
\frac {\mu}{\nu}\frac {dg}{dx} =
\sqrt {\frac {\lambda^3}{\nu}}\frac {f^{3/2}}{\sqrt {x}}.
\end{array}
$$

To hold the initial shape of the system (\ref{34}) it is necessary that
$\lambda = \mu \nu, \; \mu^2 = \lambda^3 \nu$ or, equivalentely,
$\lambda = \nu^{-3}, \; \mu = \nu^{-4}$. Hence, the system (\ref{34})
turns into itself if
$f \mapsto \nu^{-3} f, \; g \mapsto \nu^{-4} g, \; x \mapsto \nu x$.
It allows to reduce the dimension of the system. To do this, one needs
to take as new variables such combinations of old ones that are
invariant under the action of $G$. The simplest such combinations are
$y = x^3 f, \; z = x^4 g.$ In order that the system in the new
coordinates becames autonomous, (i.e. does not content the independent
variable) we put $x = \exp \rho$. In this case shifts of $\rho$
transfers into multiplications by a constants $\nu$. Then the system
(\ref{34}) takes the form

\beq
\begin{array}{l}
\frac {dy}{d\rho} = 3y + z \\
\frac {dz}{d\rho} = 4z + \sqrt {y^3}.
\end{array}
\label{35}
\eeq

The theorem is proved.

$\Box$

Now we can explore the phase portrait of the system (\ref{35}) on the
two-dimensional phase plane $(y,z)$. Besides the obvious equilibrium
at the origin we have the singular point
${\cal{A}}=(y=144, z=3\cdot 144).$ By returning to the initial
variables, the singular point $\cal{A}$ gives the famous partial
solution to the equation (\ref{31}) guessed by Sommerfeld,
$f = 144x^{-3}$. It is agreed (cf. \cite{Bet}) that this solution
gives adequate asymptotics of the exact solution for atoms as
$r \to + \infty$, but does not satisfy the boundary conditions at
the origin. This fact is understandable if to use the system (\ref{35}).
It is easy to check that the origin is the unstable knot of
the system (\ref{35}), and the point $\cal{A}$ is the saddle. The exact
solution for the atom that satisfies both boundary conditions is the
separatrix passing from the origin and tending to $\cal{A}$ as
$\rho \to + \infty$. It is not hard to show the existence of the
separatrix. And for it $r \to 0$ as $\rho \to - \infty$.

However we are interesting in boundary conditions (\ref{33}). In terms
of the system (\ref{35}) these boundary conditions transfer into the
following ones (let us recall that $y = x^3 f, \; z = x^4 g.$)

$$
y(\rho_0) = \frac {2^7 R^3}{(3\pi)^2}, \quad
y(\rho_1) = \frac {2^7 (R+r_0)^4}{Qe^2(3\pi)^2} \zeta,
$$
where

$$
\rho_0 = \ln \frac {2^{7/3}R}{(3\pi)^{2/3}}, \quad
\rho_1 = \ln \frac {2^{7/3}(R+r_0)}{(3\pi)^{2/3}}.
$$

Finally we have

$$
\begin{array}{l}
y(\ln \frac {2^{7/3}R}{(3\pi)^{2/3}}) = \frac {2^7 R^3}{(3\pi)^2}, \\
y(\ln \frac {2^{7/3}(R+r_0}{(3\pi)^{2/3}}) =
\frac {2^7 (R+r_0)^4}{Qe^2(3\pi)^2} \zeta
\end{array}
$$

Thus we need a trajectory of the phase plane $(y,z)$ that starts with
the vertical line $y=\frac {2^7 R^3}{(3\pi)^2}$ and leds to the
other vertical line
$y = \frac {2^7 (R+r_0)^4}{Qe^2(3\pi)^2} \zeta$ during the time
$\rho_1 - \rho_0 = \ln r_0 + \frac {1}{3}\ln \frac {128}{(3\pi)^2}$.
Coordinates of the final point $(y(\rho_1),z(\rho_1))$ of that
trajectory give the derivative of the potential
$\frac {\partial V}{\partial r} (R+r_0)$, and that in its turn is equal
to the shielded force of attraction of the kernel. Because of
$z=gx^4, \; g=\frac {\partial f}{\partial x}, \; f=-\frac {rV}{Qe^2}$,
we obtain

\beq
\begin{array}{c}
z(\rho_1) = - \frac {(R+r_0)^4}{b^4}\frac {\partial (rV)}{\partial r}
\frac {\partial r}{\partial x} =       \\
-\frac {(3\pi)^2}{2^7 Qe^2} (V+r\frac {\partial (V)}{\partial r}) =
\frac {(3\pi)^2}{2^7 Qe^2} (\zeta + \kappa (r_0)\frac {Qe}{(R + r_0)})
\end{array}
\label{36}
\eeq

Substitute into the equation (\ref{36}) the value $\kappa (r_0)$ from
the equation (\ref{21}). We obtain the relationship for finding the
unknown parameter $r_0$. Let us remark that it is natural to take for
$R$ values of fireballs that have been surveyed
$R \approx 10 cm$ (cf. \cite{Sta}).

The system (\ref{36}), (\ref{21}) can be solved by Newton's method.
It is possible to take for initial data nonshielded potential
$V_1(r) = - Qe^2/(R+r); \; \kappa_1 (r_0) = 1$, and sufficiently small
value
$r_0$, say, $(r_0)_1 = R/100$. Find the corresponding value
$(\rho_0), \, (\rho_1), \; y(\rho_0), \, y(\rho_1)$. Using
(\ref{35}) and reversing the time current, we obtain the trajectory
$y_1(\cdot), z_1(\cdot)$. To calculate the derivative of the first
approximation, we have to consider the variational equation of the
system (\ref{35}) for the trajectory $y_1(\cdot), z_1(\cdot)$. We
obtain the derivative of
$\kappa_1$ at the point $\kappa_1 (r_0)$. Using Newton's algorithm,
find $\kappa_2$ and new $(r_0)_2$. After that the process is repeated.

\bigskip
{\bf 4. The field of electrons.}
\bigskip

Our aim is to quantize the field of electrons itself that moves at
trajectories of the geodesic flow on the sphere. To simplify the
construction we will take the tangent bundle of unit tangent vectors
$T^1S^2$ to the two-dimensional sphere as the phase space. This gives
us the possibility disregard the electromagnetic field of the kernel,
since its attraction was properly accounted. Objects of quantization
are pairs of electrons, i.e. bosons having integral-valued spin.

One of the way is to quantize the manifold $T^1S^2$, using the
technique of geometrical quantization in terms of
Kirilov-Kostant-Souriau \cite{Ber}, \cite{Kir}, \cite{Kos},
\cite{Kum}, \cite{Raw}. \cite{Sou1}, \cite{Sou2}, \cite{Tju},
\cite{Hit}. It would be very interesting to try this approach and to
develop the technique of geometrical quantization, considering
bundles over the Kepler manifold and the projectivization of its
sections as states.

But we will use another approach. As far as electrons are moving at
trajectories of the geodesic flow, we take as the main object of the
theory (quasi-particles) geodesics itself. Consider the manifold of
nonoriented geodesics on the sphere $S^2$ \cite{Bes}. This manifold is
$\R\P^2$. Let us recall that the tangent bundle of unit tangent
vectors $T^1S^2$ to the two-dimensional sphere is $\R\P^3$.
The stereographic projection of $\R\P^3$ from the center of hemisphere
induces on $\R^3$ the metric of the elliptic space $Ell(3)$. The
geodesic flow gives by this projection the set of parallel (in the
sense of Clifford) lines. This set constitute the manifold $\R\P^2$.
The distant between Clifford-parallels remains constant in the metric
of $Ell(3)$. It is natural to take it as a distant between geodesics on
$S^2$. So, our quasi-particles are points of $\R\P^2$.

Let us find interrelation of quasi-particles. Denote the corresponding
potential by $V(\rho (x,y))$, where $\rho (x,y)$ is the distant from
$x$ to $y$ in $\R\P^2$. To be uniquelly defined on $\R\P^2 \times \R\P^2$,
the function $V(\rho)$ must be zero at infinity, i.e. $V(\pi /2) = 0$,
and besides, $V(\rho) \to \infty$ as $\rho \to 0$. This defines a
dynamical system on a set of points of $\R\P^2$. For finite sets
of points it is the system of ordinary differential equations.
Equidistant systems of points on $\R\P^2$ are stable equilibriums.
These equilibriums are not isolated, since the set of points can move
along $\R\P^2$ as a rigid body. If equidistant systems do not exist
then we have no equilibriums. But, in view of condition
$V(\rho) \to \infty$ as $\rho \to 0$, for any fixed value of energy,
distants between points during the motion admits lower estimation by a
positive constant.

Thus, the phase space is the projective plane $\R\P^2$ (recall that
spherons are parametrized by points of $\R\P^2$ also). To find the law
of interaction of quasi-particles, let us calculate the interaction of
two uniformly charged geodesics on $S^2$.

\bigskip

\begin{theorem}
The potential of interaction of two uniformly charged geodesics on
$S^2$ equals
$$
\frac {1}{\sin \theta},
$$
where $\theta$ is the angle between these geodesics.
\end{theorem}

\Proof.

First of all, find the interaction of a point $A \in S^2$ with the
uniformly charged geodesic $\sigma_1 \in S^2$. Introduce in $S^2$
the spherical coordinates
$x = \cos \phi \cos \psi, y = \sin \phi \cos \psi, z = \sin \psi.$
Without loss of generality, take $\sigma_1$ as

$$
x = \cos t, \; y = \sin t, \; z = 0, \; (0 \le t \le 2\pi).
$$

Take a point $A$ with the coordinates $(\phi, \psi)$. The vector $l$
that joins $A$ with the moving point $t$ of $\sigma_1$, is

$$
\begin{array}{l}

l = (\cos \phi \cos \psi -\cos t, \sin \phi \cos \psi -\sin t,
\sin \psi), \\

|l| = \sqrt {2 - 2\cos t \cos \phi \cos \psi -
2\sin t \sin \phi \cos \psi}. \\

\end{array}
$$

The resultant force for all points of the $\sigma_1$ is aligned with
the meridian passing through the point $A$ that is orthogonal to
$\sigma_1$, i.e. along the tangent to the circle
$(x = \cos \phi \cos s, \; y = \sin \phi \cos s, \; z = \sin s$
at the point $s = \psi.$ This is the vector
$e = (-\cos \phi \sin \psi, -\sin \phi \sin \psi, \cos \psi).$
The angle between $l$ and $e$ will be denoted by $\gamma.$ Then

$$
\begin{array}{c}

\cos \gamma = (l,e) = \frac {\cos t \cos \phi \sin \psi +
\sin t \sin \phi \sin \psi}{\sqrt {2 - 2\cos t \cos \phi \cos \psi -
2\sin t \sin \phi \cos \psi}} = \\

\frac {\sin \psi \cos (t - \phi)}
{\sqrt {2 - 2 \cos \psi \cos (t - \phi)}}. \\

\end{array}
$$

The element of the force acting on $A$ from the portion $dt$ of the
circle $\sigma_1$, equals

$$
dF(t) = \frac {l dt}{|l|^3} = \frac {l dt}
{(2 - 2 \cos \psi \cos (t - \phi))^{3/2}}.
$$

The projection on the direction of the resultant is

$$
(e,dF(t)) = \frac {\sin \psi \cos (t - \phi)}
{(2 - 2 \cos \psi \cos (t - \phi))^2}.
$$

The force acting to the point $A$ from all circle $\sigma_1$ equals

$$
\begin{array}{l}

\int_0^{2\pi} \frac {\sin \psi \cos (t - \phi) dt}
{(2 - 2 \cos \psi \cos (t - \phi))^2} = \\

2\int_0^{\pi} \frac {\sin \psi \cos (\tau ) d\tau}
{(2 - 2 \cos \psi \cos (\tau ))^2}. \\
\end{array}
$$

Rearrange this integral by the change of variable
$\tan \frac {\tau }{2} = u$
$$
\begin{array}{l}
\frac {\sin \frac {\psi }{2}}{2\cos^3 \frac {\psi }{2}}
\int_0^{\infty} \frac {(1-u^2)du}{(u^2 + \tan^2 \frac {\psi }{2})^2} = \\

\frac {\pi \sin \frac {\psi }{2}}{8\cos^3 \frac {\psi }{2}}
\left \{
\frac {1}{\tan^3 \frac {\psi }{2}} -
\frac {1}{\tan \frac {\psi }{2}}
\right \}
= \frac {\pi \cos \psi }{2 \sin^2 \psi}.

\end{array}
$$

Find the interaction of the geodesic $\sigma_1$ with that $\sigma_2$
obtained from it by rotation round the axis $Ox$ by an angle $\theta$.

$$
\sigma_2 = (\cos \tau, \sin \tau \cos \theta , \sin \tau \sin \theta).
$$

The spherical coordinates of a moving point $\tau \in \sigma_2$ are

$$
\begin{array}{l}
\tan \phi_{\tau} = \cos \theta \tan \tau; \qquad \cos \phi_{\tau} =
\frac {1}{\sqrt {1 + \cos^2 \theta \tan^2 \tau}}; \qquad
\sin \phi_{\tau} = \frac {\cos \theta \tan \tau}{\sqrt {1 +
\cos^2 \theta \tan^2 \tau}}; \\

\sin \psi_{\tau} = \sin \tau \sin \theta; \qquad
\cos \psi_{\tau} = \sqrt {1- \sin^2 \tau \sin^2 \theta}.
\end{array}
$$

Let us remark that $ \sqrt {1- \sin^2 \tau \sin^2 \theta} =
\sqrt {\cos^2 \tau + \sin^2 \tau \cos^2 \theta}$.
The force $F_{\tau}$ with the
absolute value

$$
|F_{\tau}| = \frac {\pi \cos \psi_{\tau}}{2\sin^2 \psi_{\tau}} =
\frac {\pi \sqrt {\cos^2 \tau + \sin^2 \tau \cos^2 \theta}}
{2\sin^2 \tau \sin^2 \theta}.
$$
acts (along the meridian) on a point $\tau$
We need its moment relative to the axis $O_x$. Hence, we
project $F_{\tau}$ on the normal to the plane $\sigma_2$, i.e.
on the vector

$$
\eta = \left (
\begin{array}{ccc}
1 & 0 & 0 \\
0 & \cos \theta & -\sin \theta \\
0 & \sin \theta & \cos \theta
\end{array}
\right )
\left (
\begin{array}{c}
0 \\
0 \\
1
\end{array}
\right ) =
\left (
\begin{array}{c}
0 \\
-\sin \theta \\
\cos \theta
\end{array}
\right ).
$$

The unit vector $\xi$ tangent to the meridian passing through the point
$\tau$, equals

$$
\begin{array}{c}
\xi = (-\cos \phi_{\tau} \sin \psi_{\tau}, \;
-\sin \phi_{\tau} \sin \psi_{\tau}, \; \cos \psi_{\tau}) = \\
\left (
-\frac {\sin \tau \sin \theta \cos \tau}
{\sqrt {\cos^2 \tau + \sin^2 \tau \cos^2 \theta}}, \;
-\frac {\sin \tau \sin \theta \cos \theta}
{\sqrt {\cos^2 \tau + \sin^2 \tau \cos^2 \theta}}, \;
{\sqrt {\cos^2 \tau + \sin^2 \tau \cos^2 \theta}}
\right ).
\end{array}
$$

Further on

$$
\begin{array}{c}
(\xi,\eta) = -\frac {\sin^2 \tau \sin^2 \theta \cos \tau}
{\sqrt {\cos^2 \tau + \sin^2 \tau \cos^2 \theta}} +
\cos \tau \sqrt {\cos^2 \tau + \sin^2 \tau \cos^2 \theta} = \\

\frac {\cos \theta}
{\sqrt {\cos^2 \tau + \sin^2 \tau \cos^2 \theta}}.
\end{array}
$$

This means that we have the projection

$$
|F_{\tau}|(\xi,\eta) = \frac {\pi \cos \theta}{\sin^2 \tau \sin^2 \theta}.
$$

The arm of this force is the distant from the point
$(\cos \tau, \sin \tau \cos \theta, \sin \tau \sin \theta)$ to the axis
$O_x$. It equals $\sin^2 \tau$. The moment acting on the portion
$d\tau$ of the circle $\sigma_2$ equals

$$
M_{\tau} d\tau = \frac {\pi \cos \theta }{2\sin^2 \theta}d\tau
$$

The total moment acting on the $\sigma_2$ is

$$
\int_0^{2\pi} M_{\tau} d\tau = \frac {\pi^2 \cos \theta}{\sin^2 \theta}
$$

Thus, points of the projective plane $\R\P^2$ are repelling (up to a
multiplicative constant) with the force

$$
\frac {\cos \theta}{\sin^2 \theta},
$$
where $\theta$ is the angular distant between them. The potential of
this force is

$$
\frac {1}{\sin \theta}.
$$

The theorem is proved

$\Box$
\bigskip

It is reasonable to believe that the kinetic energy of quasi-particles,
as well as for spherons, is described by the Laplace-Beltrami operator.
In the homogenous coordinats $(x,y,z)$ on the projective space, this
operator has the form

\beq
\begin{array}{l}
D = (y^2 + z^2)\frac {\partial^2}{\partial x^2} +
(z^2 + x^2)\frac {\partial^2}{\partial y^2} +
(x^2 + y^2)\frac {\partial^2}{\partial z^2} -  \\
2xy\frac {\partial^2}{\partial x \partial y} -
2yz\frac {\partial^2}{\partial y \partial z} -
2zx\frac {\partial^2}{\partial z \partial x} -   \\
2x\frac {\partial}{\partial x} -
2y\frac {\partial}{\partial y} -
2z\frac {\partial}{\partial z}.
\end{array}
\label{37}
\eeq

The distant between points of $\R\P^2$ is defined by the logarithm of
the cross ratio. In view of the preceeding theorem, the potential of
the pairwise interrelation of points with the coordinates
$(x_i,y_i,z_i)$ ¨  $(x_j,y_j,z_j)$ is the sine in power $-1$ of the
corresponding angle

$$
V_{ij} = g \frac{\sqrt{(x_i^2+y_i^2+z_i^2)(x_j^2+y_j^2+z_j^2)}}
{\sqrt{(x_iy_j-x_jy_i)^2 + (y_iz_j-y_jz_i)^2 + (z_ix_j-z_jx_i)^2}},
$$
where $g$ is the interaction parameter.

Denote by $D_i$ the Laplace-Beltrami operator in $Ell(\R^2)$ applying
to the coordinates $q_i$. The Hamiltonian of the system that consists of
$N$ interacting quasi-particles is

\beq
H^1 = - \gamma \sum_{i=1}^N \frac {1}{2}D_i + \sum_{i<j}V{ij}.
\label{38}
\eeq

As it was expected, the Hamiltonian is homogeneous. We consider it on
the space $\R\P^2$, i.e. on functions of zero degree homogeneity
relative to $x_i,y_i,z_i$ for any $i$. Besides, it is invariant
relative to the group of permutations of points. We obtain an analogue
of the Calogero-Sutherland-Moser problem (CSM) \cite{Sut}. But here
the part of the circle plays the projective space $\R\P^2$. It would be
good to find the asymptotics of the spectrum as $N \to \infty$.

We will use the Euler formulas of the first and second orders for
functions of zero degree homogeneity.

$$
\begin{array}{l}
x\frac {\partial f}{\partial x} +
y\frac {\partial f}{\partial y} +
z\frac {\partial f}{\partial z} = 0; \\
x^2\frac {\partial^2 f}{\partial x^2} +
y^2\frac {\partial^2 f}{\partial y^2} +
z^2\frac {\partial^2 f}{\partial z^2} +
xy\frac {\partial^2 f}{\partial x \partial y} +
yz\frac {\partial^2 f}{\partial y \partial z} +
zx\frac {\partial^2 f}{\partial z \partial x} = 0.
\end{array}
$$

By adding to and subtracting from $D$ the value

$$
x^2\frac {\partial^2}{\partial x^2} +
y^2\frac {\partial^2}{\partial y^2} +
z^2\frac {\partial^2}{\partial z^2},
$$
and then by using the Euler formulas, we reduce (\ref{37}) to the form

$$
D_i =  |q_i|^2
\left (
\frac {\partial^2}{\partial x_i^2} +
\frac {\partial^2}{\partial y_i^2} +
\frac {\partial^2}{\partial z_i^2}.
\right )
$$

Hence, the operator $H^1$ takes the form

\beq
H^1 = - \frac {\gamma }{2} \sum_{i=1}^N |q_i|^2 \Delta_i  + \sum_{i<j}V{ij}.
\label{39}
\eeq

Enumerate eigenvalues of the Hamiltonian (\ref{39}) by natural numbers
$k$ and denote them by $-\sigma_k$. The corresponding complete
orthonormal system of eigenfunctions is denoted by
$Y_{k}(S),\; S \in (S^2)^N$.
Expand a solution of the equation

\beq
\frac {\partial^2 U}{\partial t^2} = H^1(U)
\label{40}
\eeq
in term of functions $Y_k(S)$. Coefficients of the series are defined by
expression

\beq
\hat U_k(t) = \int_{(S^2)^N} U(S,t)Y_k(S)\, dS.
\label{41}
\eeq

Apply the operator $\frac {\partial^2}{\partial t^2}$ to (\ref{41}).
Since the operator $H^1$ is selfajoint and $Y_k(S)$ are its
eigenfunctions with the eigenvalues $- \sigma_k$, we obtain

$$
\frac {\partial^2 \hat U_k(t)}{\partial t^2} =
\int_{(S^2)^N} H^1 U(S,t) \, Y_k(S) dS =
\int_{(S^2)^N} U(S,t) \, H^1 Y_k(S) dS = - \sigma_k\hat U_k(t).
$$

Hence, we get the disjoin system of ordinary differential equations

\beq
\frac {\partial^2 \hat U_k(t)}{\partial t^2} = -
\sigma_k\hat U_k(t).
\label{42}
\eeq

Each equation of the system (\ref{42}) is the equation of the harmonic
oscillator with the unit mass and the frequency $\sqrt{\sigma_k}$
(we shall call it by the oscillator $k$). The application of the
standard quantization proceedure leads to the secondary quantization of
the quasi-partical oscillations. The Hamiltonian of the oscillator $k$
is

$$
H^1_k = \frac {1}{2} (p^2 + \sigma_kq^2).
$$

We anumerate eigenfunctions of the operator $H^1_k$ by indices $n$ and
denote by $\Xi_n (k)$. It is known that the corresponding
eigenvalues are

$$
E_n (k)= (n + \frac {1}{2})\sqrt{\sigma_k}.
$$

We introduce operators as in (\ref{17}).

\beq
A^+(k) \Xi_n (k) = \sqrt{n + 1} \, \Xi_{n + 1}(k), \quad
A^-(k) \Xi_n (k) = \sqrt{n} \, \Xi_{n - 1}(k).
\label{43}
\eeq

The ground-state of the system corresponds to $n = 0$ since

$$
A^-(k) \Xi_0(k) = 0,
$$
and

\beq
\Xi_n (k) = \frac {(A^+(k))^n}{\sqrt{n !}} \Xi_0(k).
\label{44}
\eeq

The formulas (\ref{43}) show that states $\Xi_n (k)$ are
eigenfunctions of the operators $A^-(k)A^+(k)$ and $A^+(k)A^-(k)$
with the eigenvalues $n + 1$ and $n$
correspondingly. Consequently, it is natural to call the operator
$A^+(k)A^-(k)$ by the operator of number of quasi-particals
of type $k$. Denote this operator by $\hat n^q (k).$ The
commutation relations take the standard form

$$
[A^-(k),A^+(k)] = 1, \; [A^-(k),A^-(k)] =
[A^+(k),A^+(k)] = 0.
$$

Meanwhile operators that relate to noncoinciding values of indices
commute. The Hamiltonian for an oscillator is

$$
H^1(k) = \sqrt{\sigma_k}(A^+(k)A^-(k) + \frac {1}{2}{\rm I}^q(k)),
$$
or
$$
H^1(k) = \sqrt{\sigma_k}(\hat n^q + \frac {1}{2}{\rm I}^q(k)),
$$
where ${\rm I}^q(k)$ is the identity operator on the space of states
of the oscillator $k$. The full Hamiltonian of the field of the
quasi-papticles is represented by the sum

\beq
H^1 = \sum_k H^1(k) =
\sum_k \sqrt{\sigma_k}(\hat n^q + \frac {1}{2}{\rm I}^q(k)).
\label{45}
\eeq

\bigskip

Let us find the operator $H^{int}$ of interaction of spherons with
quasi-particles. The interaction annihilates a
quasi-particle of an energy $E_3$ that corresponds to the oscillator
$k_3$, and this excites a wave of an energy $E_1$ that corresponds to
the spheron $(m_1,2l_1)$. This wave create quasi-particle
of an energy $E_2$ that corresponds to the oscillator $k_2$, and this
annihilates a wave of an energy $E_4$ that corresponds to the spheron
$(m_4,2l_4)$. So, we can write the operator of interaction $H^{int}$.

\beq
\begin{array}{c}
H^{int} = \\ - W \sum_{k_j,l_i,|m_i| \le 2l_i}
(a^+(m_1,2l_1)A^+(k_2)A^-(k_3)a^-(m_4,2l_4))
\delta(E_1+E_2-E_3-E_4),
\end{array}
\label{46}
\eeq
where $E_i$ are energies of the corresponding spherons and
quasi-particles. The $\delta$-function at the right hand side of
(\ref{46}) provides energy conservation.

The full Hamiltonian is

\beq
\begin{array}{c}
H = H^0 + H^1 + H^{int} = \\
\sum_{l,|m| \le 2l}\sqrt{2l(2l+1)}(\hat n^s (m,2l) +
\frac {1}{2}{\rm I}^s(m,2l)) +
\sqrt{\sigma_k}(\hat n^q + \frac {1}{2}{\rm I}^q(k)) - \\
W\sum_{k_j,l_i,|m_i| \le 2l_i}
(a^+(m_1,2l_1)A^+(k_2)A^-(k_3)a^-(m_4,2l_4))
\delta(E_1+E_2-E_3-E_4).
\end{array}
\label{47}
\eeq

\bigskip

The problem is significantly simplified if we neglect terms of
Coulomb interaction $V_{ij}$ considering free
quasi-particles as a first approximation.

At first, take the special case when $\gamma = 1$. In this case
both Hamiltonians $H^0$ and $H^1$ act on isomorphic spaces that we
denote by ${\goth A}^s$ and ${\goth A}^q$ correspondingly. Both these
operators are diagonal in basis consisted of spherical harmonics.
Its spectra coincide with N-multiple spectrum of the harmonic
oscillator. Note that operators $a^{\pm}$ and $A^{\pm}$ commute,
acting on different objects (different components of a
state-vector) The Hamiltonian equals

\beq
\begin{array}{c}
H = \sum_{l,|m| \le 2l} \sqrt{2l(2l+1)}(\hat n^s (m,2l) +
\frac {1}{2}{\rm I}^s(m,2l)) + \hat n^q (m,2l) +
\frac {1}{2}{\rm I}^q(m,2l)) -  \\
W\sum_{l_i,|m_i| \le 2l_i}
(a^+(m_1,2l_1)A^+(m_2,2l_2)A^-(m_3,2l_3)a^-(m_4,2l_4))
\delta(l_1+l_2-l_3-l_4).
\end{array}
\label{48}
\eeq

We shall call the interaction switch-back if
$m_1=m_3,\, l_1=l_3$ and $m_2=m_4,\, l_2=l_4$. Find the spectrum of
Hamiltonian under the natural conjecture of switch-back interaction.
In this case the interaction behaves as if it were two different acts.
The first act is the annihilation of a quasi-partical of type
$(m_1,2l_1)$ with the synchronous excitation of a spheron of the same
type. The second one is the annihilation of a spheron of type
$(m_2,2l_2)$ with the birth of a quasi-partical of the same type.
Since the Hamiltonian contains the sum of all these processes, one can
imagine that the state-vector of quasi-particles turns into the
state-vector of spherons and conversely. The interaction
Hamiltonian is

\beq
\begin{array}{c}
H = \sum_{l,|m| \le 2l} \sqrt{2l(2l+1)}(\hat n^s (m,2l) +
\frac {1}{2}{\rm I}^s(m,2l)) + \\
\sum_{l,|m| \le 2l} \sqrt{2l(2l+1)}
(\hat n^q (m,2l) + \frac {1}{2}{\rm I}^q(m,2l)) -  \\
W\sum_{l_i,|m_i| \le 2l_i}
(a^+(m_1,2l_1)A^+(m_2,2l_2)A^-(m_1,2l_1)a^-(m_2,2l_2)).
\end{array}
\label{49}
\eeq

Let us rewrite the operator $H$ given by the formula (\ref{49}), by
using the polarization of the state-space
$\goth A = {\goth A}^s \oplus {\goth A}^q$. Take in both isomorphic
subspaces ${\goth A}^s$ and ${\goth A}^q$ the concordant bases. Then
the operator $H^{int}$ interchanges vectors of these subspaces and
multiply that by $-W$, i.e.
$H^{int}: (x,y) \mapsto (-Wy,-Wx)$. Consequently, the operator $H$ has
the block form:

\beq
\left (
\begin{array}{cc}
B_{11} & B_{12} \\
B_{21} & B_{22}
\end{array}
\right ),
\label{51}
\eeq
where $B_{11}=B_{22}=Diag (b_1,b_2,\dots), \; B_{12}=B_{21}=-WI$.

\bigskip

To find the spectrum of the operator $H$, we explore a little more
general problem. Consider an operator $A$, acting on the direct sum of
two isomorphic separable Hilbert spaces $K \oplus K$, with the block
structure

\beq
A =
\left (
\begin{array}{cc}
B_1 & B_2 \\
B_2 & B_1
\end{array}
\right ),
\label{52}
\eeq
where $B_i$ are commuting selfajoint operators acting on the space $K$.
In view of the von-Neumann theorem \cite{Neu}, there exists selfajoint
operator $C$ such that both operators $B_i$ are functions of $C$. Take
$E_{\lambda}$ --- the spectral decomposition of the space $K$,
corresponding to the operator $C$. Then

$$
B_1 = \int \alpha (\lambda) dE_{\lambda}, \quad
B_2 = \int \beta (\lambda) dE_{\lambda}.
$$

Note that the diagonal subspaces
$K_1 = \{(x,x)\}$ and $K_2 = \{(x,-x)\}$ are invariant relative to the
operator $A$. By using the polarization $K \oplus K = K_1 \oplus K_2$,
we reduce the operator $A$ to the block-diagonal form

\beq
A =
\left (
\begin{array}{cc}
(B_1+B_2)  &     0     \\
    0      & (B_1-B_2)
\end{array}
\right ).
\label{53}
\eeq

Denote by $E_{\lambda} \oplus E_{\lambda}$ the operators that
correspond to the spectral family $E_{\lambda}$ into the subspace
$K_1 = \{(x,x)\}$, and by
$E_{\lambda} \ominus E_{\lambda}$ the operators that
correspond to the spectral family $E_{\lambda}$ into the subspace
$K_2 = \{(x,-x)\}$.

We proved the following theorem

\bigskip

\begin{theorem}
The spectral representation of the operator $A$ on the space $K_1$ is

$$
A|_{K_1} = \int (\alpha (\lambda) + \beta (\lambda))
d(E_{\lambda} \oplus E_{\lambda}).
$$

The spectral representation of the operator $A$ on the space $K_2$ is

$$
A|_{K_2} = \int (\alpha (\lambda) - \beta (\lambda))
d(E_{\lambda} \ominus E_{\lambda}).
$$
\end{theorem}

\bigskip

\begin{corollary}
The spectrum of the operator $H$ can be obtained from the
spectrum of the operator $B_{11}=B_{22}$ by the shifts (as a whole)
both to the right and to the left for a distance of $W$.
\end{corollary}

\bigskip

Now consider the case when $\gamma = n/m$ any rational number. As
before we neglect the interactions defined by terms $V_{ij}$. Let us
unify $n$ exciting states of spherons in a cluster that we
call pseudo-spheron. In the same way, $m$ exciting states of
quasi-particles will be a cluster called pseudo-particle. For
sufficiently many quasi-particles and quasi-spherons, the Hamiltonians
$H^0$ and $H^1$ will asymptotically coincide. Irrational values of
$\gamma $ can be approximated by rational ones preserving the
asymptotics of the spectrum. So, the problem can be reduced to the
special case when $\gamma = 1$.

\bigskip

Since the spectrum of the operator $B_{11}=B_{22}$ is known, then at
the approximation in question we obtain the spectrum of the
Hamiltonian $H$. So, we obtain a possibility (in principle) to find
mean values of important operators. Besides, our model Hamiltonian
provides the basis for application of the perturbation theory for more
realistic Hamiltonians.

\bigskip
\bigskip

      The question:
\bigskip

Does there exist a coherence of electrons under a line (usual)
lightning? To resolute the question on the quantum character of the
usual lightning (or any spark discharge) one needs to explore
oscillations of cylindrical flow of electrons. Propose some
conjectures concerning this question.

\bigskip

{\bf 5. Conjectures}

\bigskip

The inner canal of a lightning is the ionized, positively charged
plasma. Electrons are winding round this canal into a
cylindrical cover. Stepped leader is, apparently, a process of
transformation into a superconductive state. And the
superconductive state itself is connected with macroscopic separation
of positive and negative plasma. It seems that Abrikosov's threads in
the theory of usual (low-temperature) superconductivity have the same
structure: canal from ions surrounded by a cover of winding electrons
(not for nothing these threads are called vortical). It is precisely
this separation of plasma which gives rise to superconductivity. The ionized
air in the course of thunderstorms behaves like a superconductor of
the second kind. Stepped leader is the process of generation of analogs
of Abrikosov's threads --- media of plasma superconductivity. These
threads can curve, branch, and terminate due to rushes of the wind.
The lightning flying through these canals repeates its form. The
stepped leader constructs the canal of lightning by steps, as if it
makes links of a chain. So, one of the possible ways of forming
fireballs is separation of one of the chain from the canal of
lightning. Nice recent theory of lightning (Runaway Breakdown)
\cite{Gur} takes into account the leading role of high-energy cosmic
rays. It seems that Runaway Breakdown relates rather to the stepped
leader then to the lightning itself. This is in agreement with the
recent observations \cite{Dwy}, where it was shown that steps of the
stepped leader were well correlated with the high-energetic
microsecond burst of gamma emission.

The question arises: why the lightning extincts so fast if there
exists the canal of superconductivity? Perhaps it happens due to the
magnetic field. Too strong magnetic field, that is excited by the
powerfull current in the superconductor, destruct its creator --- the
state of superconductivity that induces this current. Electrons begin
to fall from the cover into the canal, superconductivity disappears,
and the lightning extincts. The difference in pressure generates the
shock wave --- the thunder.

\end{document}